# From harmonic mappings to Ricci solitons

Sergey Stepanov[1,2], Irina Tsyganok[2]

[1] Department of Mathematics, All Russian Institute for Scientific and Technical Information of the Russian Academy of Sciences, 20, Usievicha street, 125190 Moscow, Russia,
e-mail: s.e.stepanov@mail.ru

[2] Department of Mathematics, Finance University under the Government of the Russian Federation, Leningradsky Prospect, 49-55, 125468 Moscow, Russia,
e-mail: i.i.tsyganok@mail.ru

**Abstract.** The paper is devoted to the study of the global geometries of harmonic mappings and infinitesimal harmonic transformations and presents their applications to the theory of Ricci solitons.

## 1. Introduction

In this paper we study the global geometry of harmonic mappings ([1]; [2]; [3]) and, in particular, infinitesimal harmonic transformations of Riemannian manifolds ([4] [5]) and present their applications to the Ricci soliton theory ([6] [7]).

Our results will be obtained using *the Bochner technique* and, in particular, the maximum principles of Hopf ([8, Theorem 1]; [9, Theorem 2.1) and Bochner ([9, Theorem 2.2]). We will also use Yau, Li and Schoen results on the connections between the geometry of a complete smooth manifold and the global behavior of its subharmonic functions (see, for example, [10] and [11]).

In the paper we continue the study of certain connections between the theory of infinitesimal harmonic transformations and the theory of Ricci solitons which we began in [12]. Theorems and corollaries of the paper complement our results from [5]; [12]-[16].

In the second section of the paper, we give brief survey of basic facts of the theory of harmonic mappings between Riemannian manifolds in the large. Results of the third section of our paper with the title "Infinitesimal harmonic transformations" are obtained as analogs of results of the second section of the paper. In turn, the results of the fourth section which has the title "Ricci solitons" are applications of the results of the third section of our paper.

The present paper is based on our report at the conference "Differential Geometry" organized by the Banach Center in June 18-24, 2017 in Bedlewo (Poland) and on our plenary lecture at the International Conference "Modern geometry and its applications" organized by the Lobachevsky Institute of Mathematics and Mechanics from November 27 to December 3, 2017 in Kazan (Russia).

## 2. Harmonic mappings in the large

In this section we shall give brief survey of basic facts of the differential geometry of harmonic mappings between Riemannian manifolds in the large (see, for example, [1]; [2] and [3]). Moreover, we shall have given priority to new proofs of the well-known vanishing theorems of harmonic mappings.

Let $(M, g)$ and $(\overline{M}, \overline{g})$ denote complete Riemannian manifolds of dimensions $n$ and $\overline{n}$, respectively. We call the *energy density* of a smooth map $f : (M, g) \to (\overline{M}, \overline{g})$ the nonnegative scalar function $e(f) : M \to \mathbb{R}$ such that $e(f) = \frac{1}{2} \| f_* \|^2$ where $\| f_* \|^2$ denotes the squared norm of the differential $f_*$, with respect to the induced metric $\widetilde{g}$ on the vector bundle $T^*M \otimes f^*T\overline{M}$ by $g$ and $\overline{g}$ (see also [1]).

It is well known that $f : (M, g) \to (\overline{M}, \overline{g})$ is a *harmonic mapping* if and only if it satisfies the *Euler-Lagrange equation*

(2.1) $$\operatorname{trace}_g \left( \widetilde{\nabla} f_* \right) = 0$$

where $\widetilde{\nabla} = \nabla \oplus \overline{\nabla}$ is the canonical connection in the vector bundle $T^*M \otimes f^*T\overline{M}$ ([1, p. 117]). Moreover, if $f$ is a harmonic mapping then from (2.1) one can obtain the following equation (see also [1; p. 123]):

(2.2) $$\Delta e(f) = \left\| \widetilde{\nabla} f_* \right\|^2 + Q(f)$$

for the *Laplace–Beltrami operator* $\Delta = \operatorname{div} \nabla$ and the scalar function

(2.3) $$Q(f) = g\left( \operatorname{Ric}, f^*\overline{g} \right) - \operatorname{trace}_g \left( \operatorname{trace}_g \left( f^*\overline{R} \right) \right)$$

where $\overline{R}$ is the Riemannian curvature tensor of $(\overline{M}, \overline{g})$ and $Ric$ is the Ricci tensor of $(M, g)$.

Let $U \subset M$ be a connected open domain. We suppose that $Q(f)$ is non-negative everywhere in $U$ and $Q(f)$ is positive in at least one point of $U$. In this case we call $Q(f)$ *quasi-positive* scalar function defined on $U \subset M$. From (2.3) we obtain the conditions when $Q(f)$ is quasi-positive on $U \subset M$ (see also [1; pp. 124-125]). First, $(M, g)$ must be a Riemannian manifold with *quasi-positive Ricci curvature* on $U$, i.e. the Ricci curvature is nonnegative in all directions at an arbitrary point of $U$ and, in addition, it is strictly positive in all directions at some point of $U$ ([17]). Second, $(\overline{M}, \overline{g})$ must be a Riemannian manifold with nonpositive sectional curvature $\overline{sec}$ at an arbitrary point of $f(U) \subset \overline{M}$, i.e. $\overline{g}(\overline{R}(X,Y)Y, X) \leq 0$ for all for all unite orthogonal vectors $X, Y \in T_x \overline{M}$ at an arbitrary point $x \in f(U) \subset \overline{M}$. In this case, the energy density function $e(f)$ satisfies the inequality $\Delta e(f) \geq 0$ at each point of $U$, by (2.2). Therefore, $e(f)$ is a *subharmonic function*. Now suppose that the energy density function $e(f)$ attains a local maximum value at some point $x \in U$, then $e(f)$ is a constant $C$ in $U$, by the *Hopf's maximum principle* ([8, Theorem 1]; [9, Theorem 2.1). If $C > 0$, then $grad\ f$ is nowhere zero. Now, at a point where the $Q(f)$ is positive, the left side of (2.2) is zero while the right side is positive. This contradiction shows that $C = 0$ and hence $f$ is constant in $U$. Thus we have proved the following statement.

**Lemma 2.1**. *Let $f : (M, g) \to (\overline{M}, \overline{g})$ be a harmonic mapping such that its energy density $e(f) = \frac{1}{2}\|f_*\|^2$ has local maximum at some point $x$ in a connected open domain $U \subset M$. If, in addition, either one of the two following conditions is satisfied: the sectional curvature $\overline{sec}$ of $(\overline{M}, \overline{g})$ is nonpositive at an arbitrary point of $f(U) \subset \overline{M}$ and $(M, g)$ has the quasi-positive Ricci curvature $Ric$ in $U$, then $f$ is constant in the domain $U$.*

On the other hand, we know from [18] other conditions when $Q(f)$ is quasi-positive on $U \subset M$. Namely, if the sectional curvature $\overline{sec}$ of $(\overline{M},\overline{g})$ is nonnegative at an arbitrary point of $f(U) \subset \overline{M}$ and $Ric \geq f^*\overline{Ric}$ at each point of $U$ for the Ricci tensors $\overline{Ric}$ of $(\overline{M},\overline{g})$, then $Q(f) \geq 0$. Furthermore, if there is at least one point of $U$ at which $Ric > f^*\overline{Ric}$, then $Q(f) > 0$ at this point. In this case, $Q(f)$ is also quasi-positive scalar function defined in $U$. Thus we can formulate the second statement.

**Lemma 2.2**. *Let $f:(M,g) \to (\overline{M},\overline{g})$ be a harmonic mapping such that its energy density $e(f) = \frac{1}{2}\|f_*\|^2$ has local maximum at some point $x$ in a connected open domain $U \subset M$. If, in addition, the sectional curvature $\overline{sec}$ of $(\overline{M},\overline{g})$ is nonnegative at an arbitrary point of $f(U) \subset \overline{M}$ and $(M,g)$ has the Ricci tensor $Ric$ such that $Ric \geq f^*\overline{Ric}$ at each point of $U$, and there is at least one point of $U$ at which $Ric > f^*\overline{Ric}$, then f is constant in the domain $U$.*

Now, suppose that $(M,g)$ is a compact Riemannian manifold, then there exists a point $x \in M$ at which the function $e(f)$ attains the maximum. At the same time, let $e(f)$ satisfies the condition $\Delta e(f) \geq 0$ everywhere in $(M,g)$. In this case, we can use the *Bochner maximum principle* which we deduce from the Hopf maximum principle. Namely, it is well known that an arbitrary subharmonic function is a constant function on a compact Riemannian manifold ([9, Theorem 2.2]). As a result, we can formulate the following theorem as a corollary of our Lemma 1 ([1]; [3, p. 23]).

**Theorem 2.1**. *Let $f:(M,g) \to (\overline{M},\overline{g})$ be a harmonic mapping between Riemannian manifolds $(M,g)$ and $(\overline{M},\overline{g})$. Assume that the sectional curvature $\overline{sec}$ of the second smooth manifold $(\overline{M},\overline{g})$ is nonpositive in every point of $f(M)$ and the first smooth manifold $(M, g)$ is a compact smooth manifold with the nonnegative Ricci tensor $Ric$. Then f is totally geodesic mapping. Furthermore, if there is at least one point of M in which its Ricci tensor is positive definite, then f is a constant mapping.*

On the other hand, using our Lemma 2 and the Bochner maximum principle ([9, Theorem 2.2]), we conclude that the following theorem holds ([18]).

**Theorem 2.2**. *Let $f:(M,g)\to(\overline{M},\overline{g})$ be a harmonic mapping between Riemannian manifolds $(M,g)$ and $(\overline{M},\overline{g})$. Assume that the sectional curvature $\overline{sec}$ of the second manifold $(\overline{M},\overline{g})$ is nonnegative in every point of $f(M)$ and the first manifold $(M,g)$ is a compact smooth manifold with the Ricci tensor $Ric \geq f^*\overline{Ric}$. Then $f$ is totally geodesic mapping with constant energy density $e(f)$. Furthermore, if there is at least one point of $M$ in which $Ric > f^*\overline{Ric}$, then $f$ is a constant mapping.*

In turn, Li and Schoen have proved in [11] that there is no non-constant, non-negative $L^p$-integrable $(0<p<\infty)$ subharmonic function $u$ on any complete Riemannian manifold $(M,g)$ with nonnegative Ricci tensor. In other word, if we suppose that $Ric \geq 0$ and $\int_M \|u\|^p dVol_g < \infty$ for a complete Riemannian manifold $(M,g)$, then $u = C$ for some constant $C$. In this case, we have $C^p \int_M dVol_g < \infty$. If $C > 0$, $u$ is nowhere zero and the volume of $(M,g)$ is finite. Side by side, we know from [10] that every complete non-compact Riemannian manifold $(M,g)$ with nonnegative Ricci tensor has infinite volume. This contradiction shows $C = 0$ and hence $u \equiv 0$. Therefore, we can formulate the following

**Lemma 2.3**. *Let $(M,g)$ be a complete non-compact Riemannian manifold with nonnegative Ricci tensor, then there is no nonzero non-negative $L^p$-integrable $(0<p<\infty)$ subharmonic function on $(M,g)$.*

Let $(M,g)$ be a complete non-compact Riemannian manifold. Given a smooth map $f:(M,g)\to(\overline{M},\overline{g})$ we definite its *energy* as ([13]):

$$E(f) = \int_M e(f) dVol_g .$$

The energy $E(f)$ can be infinite or finite. For example, $E(f) < +\infty$ for the compact Riemannian manifolds $(M,g)$ and $(\overline{M},\overline{g})$ ([2]). Then based on Lemma 2.3, we can formulate the following theorem for harmonic mapping of a complete Riemannian manifold (see also [2]; [3, p. 25]).

**Theorem 2.3**. *Let $(M,g)$ be a complete non-compact manifold with nonnegative Ricci curvature, and $(\overline{M},\overline{g})$ a Riemannian manifold with nonpositive sectional curvature $\overline{sec}$. If $f:(M,g)\to(\overline{M},\overline{g})$ is a harmonic mapping with finite energy, then $f$ is a constant map.*

Side by side, we can formulate an alternative theorem for harmonic maps from complete Riemannian manifolds to Riemannian manifolds with nonnegative sectional curvature (see also [18]). In this case, we can also use Lemma 2.3 on nonzero non-negative $L^p$-integrable $(0<p<\infty)$ subharmonic functions on a complete Riemannian manifold with nonnegative Ricci curvature for the proof of this theorem.

**Theorem 2.4**. *Let $f:(M,g)\to(\overline{M},\overline{g})$ be a harmonic mapping with finite energy. If $(M, g)$ is a complete non-compact Riemannian manifold with Ricci tensor $Ric \geq f^*\overline{Ric}$, where $\overline{Ric}$ denotes the Ricci tensor of $(\overline{M},\overline{g})$, and $(\overline{M},\overline{g})$ is a Riemannian manifold with nonnegative sectional curvature $\overline{sec}$ in every point of $f(M)$, then $f$ is a constant map.*

### 3. Infinitesimal harmonic transformations

The main results of this section are obtained as analogs of results of the second section of the paper.

A vector field $\xi$ on a complete Riemannian manifold $(M,g)$ is called an *infinitesimal harmonic transformation* of $(M,g)$ if $\xi$ generates a flow which is a local one-parameter group of harmonic transformations (in the other words, local harmonic diffeomorphisms) ([4]). Analytic characteristic of such vector field has the form

$$trace_g(L_\xi \nabla) = 0$$

for the Lie derivative $L_\xi$ in the direction of $\xi$. This formula is an analogue of the formula (2.1). In addition, we have proved in [5] that a vector field $\xi$ is an infinitesimal harmonic transformation if and only if

$$\tilde{\Delta}\theta = 2\,Ric(\xi,\cdot)$$

for the 1-form $\theta$ corresponding to $\xi$ under the duality defined by the metric $g$ and the *Hodge–de Rham Laplacian* $\tilde{\Delta}$ ([19, p. 158]).

In accordance with the theory of *harmonic maps* ([1]) we define the *energy density* of the flow on $(M,g)$ generated by an infinitesimal harmonic transformation $\xi$ as the scalar function $e(\xi) = \frac{1}{2}\|\xi\|^2$ where $\|\xi\|^2 = g(\xi,\xi)$. Then the Laplace–Beltrami operator $\Delta e(\xi)$ for the energy density $e(\xi)$ of an infinitesimal harmonic transformation $\xi$ has the form ([13]; [15])

(3.1) $$\Delta e(\xi) = \|\nabla \xi\|^2 - Ric(\xi,\xi).$$

The formula (3.1) is an analogue of the formula (2.2) for the energy density $e(f)$ of a harmonic map $f$. In this case, the following theorem is true.

**Theorem 3.1**. *Let $(M,g)$ be a Riemannian manifold and $U \subset M$ be a connected, open domain. If the energy density of the flow $e(\xi) = \frac{1}{2}\|\xi\|^2$ generated by an infinitesimal harmonic transformation $\xi$ has a local maximum in some point of U and the Ricci tensor Ric of $(M,g)$ is quasi-negative in $U$, then $\xi$ is identically zero everywhere in $U$.*

**Proof**. Let $Ric < 0$ everywhere in a connected, open domain $U \subset M$, then the energy density function $e(\xi)$ satisfies the inequality $e(\xi) \geq 0$, by (3.8). This means that $e(\xi)$ is a subharmonic function. Suppose now that the energy density function $e(\xi)$ attains a local maximum value in a point $x \in U$, then $e(\xi)$ is a constant $C$ in $U$, by Hopf's maximum principle ([8, Theorem 1]; [9, Theorem 2.1]). If $C > 0$, $\xi$ is nowhere zero. Now, at a point where the $Ric$ is negative, the left side of (3.1) is zero while the right side is positive. This contradiction shows $C = 0$ and hence $\xi \equiv 0$ everywhere in the domain $U$.

**Remark 3.1.** Theorem 3.1 is a direct generalization of the Theorem 4.3 presented in Kobayashi's monograph on transformation groups ([19, p. 57]) and Wu's proposition on a Killing vector field whose length achieves a local maximum ([17]).

As an analogue of Theorem 2.1, we can formulate the following theorem, which can be proved using the *Bochner maximum principal* ([9, Theorem 2.2]).

**Theorem 3.2.** *A compact Riemannian manifold $(M,g)$ with quasi-negative Ricci curvature doesn't have nonzero infinitesimal harmonic transformation.*

Next, we recall that the *kinetic energy* $E(\xi)$ of the flow on $(M,g)$ generated by a vector field $\xi$ is determined by the following equation ([25, p. 2])

$$E(\xi) = \int_M e(\xi) \, dVol_g.$$

**Remark 3.2.** The definition given above is consistent with the theory of harmonic mappings in the case of an infinitesimal harmonic transformation. Moreover, the energy $E(\xi)$ can be infinite and finite. For example, $E(\xi) < +\infty$ for a smooth complete vector field $\xi$ on a compact Riemannian manifold $(M,g)$.

As an analogue of Theorem 2.3, we formulate the following

**Theorem 3.3.** *Let $(M,g)$ be a complete Riemannian manifold $(M,g)$ with the nonpositive Ricci curvature. Then every infinitesimal harmonic transformation with finite kinetic energy is parallel. If, in addition, the volume of $(M,g)$ is infinite or the Ricci curvature is negative at some point of $M$, then the infinitesimal harmonic transformation is identically zero on $(M,g)$.*

**Proof.** For the proof we use the well known *second Kato inequality* ([26, p. 380])

$$-\|\xi\| \Delta \|\xi\| \leq g(\overline{\Delta}\theta, \theta)$$

where $\overline{\Delta} := -trace_g \nabla \circ \nabla$ is the *rough Laplacian* and $\theta$ is the 1-form corresponding to $\xi$ under the duality defined by the metric $g$. It is well known that the rough Laplacian satisfies the *Weitzenböck formula* ([3, p. 21]; [19, p. 44]; [26, p. 378])

$$\overline{\Delta}\theta = \widetilde{\Delta}\theta - S\xi$$

where $S$ is the *Ricci operator* defined by $g(SX,Y) = Ric(X,Y)$ for any tangent vector fields $X$ and $Y$. Therefore, the second Kato inequality can be rewritten in the form

(3.2) $$2\sqrt{e(\xi)}\,\Delta\sqrt{e(\xi)} \geq -g(\widetilde{\Delta}\theta,\theta) + Ric(\xi,\xi).$$

where $\|\xi\| = \sqrt{2e(\xi)}$. On the other hand, we have proved in [1] and [13] that $\xi$ is an infinitesimal harmonic transformation on $(M,g)$ if and only if $\widetilde{\Delta}\xi = 2S\xi$. Therefore, we obtain from (3.2) the following equation

(3.3) $$\sqrt{e(\xi)}\,\Delta\sqrt{e(\xi)} = -\frac{1}{2}Ric(\xi,\xi).$$

Let us require the Ricci tensor $Ric$ to be nonpositive. In [10, p. 664] and [27] it was shown that every nonnegative smooth function $u$ defined on a complete Riemannian manifold $(M,g)$ and satisfying the conditions $u\Delta u \geq 0$ and $\int_M u^p\,dVol_g < +\infty$ for all $p \neq 1$, must be constant. In particular, if the volume of $(M,g)$ is infinite and $u = const$, then $u = 0$. Therefore, if the smooth manifold $(M,g)$ is complete and

(3.4) $$E(\xi) = \int_M e(\xi)\,dVol_g < +\infty,$$

then the function $\sqrt{e(\xi)}$ is constant. At the same time, we obtain from (3.4) that the volume of $(M,g)$ is finite. Thus it follows from (3.2) that $\nabla\xi = 0$. On the other hand, if we suppose that $Ric_x < 0$ in a point $x \in M$, then this inequality contradicts the equation (3.1). The proof is complete.

**Remark 3.3**. If $\xi$ is nonzero vector field such that $Ric(\xi,\xi) \leq 0$, then (3.3) implies the inequality $\Delta\sqrt{e(\xi)} \geq 0$. This means that $\sqrt{e(\xi)}$ is a subharmonic function. Therefore, the result given above is an analogue of the result which we showed in Theorem 2.3.

The following corollary is valid.

**Corollary 3.1**. *Let $(M,g)$ be a connected complete noncompact Riemannian manifold of dimension $n \geq 2$ with irreducible holonomy group $Hol(g)$ and nonpositive Ricci curvature. Then each infinitesimal harmonic transformation on $(M,g)$ with finite kinetic energy identically vanishes.*

**Proof**. By Theorem 3.3, an infinitesimal harmonic transformation $\xi$ with $E(\xi) < +\infty$ on a connected complete noncompact Riemannian manifold $(M,g)$ with nonpositive Ricci curvature is parallel. Under the assumption that holonomy group $\text{Hol}(g)$ is irreducible, this relation means that $\xi \equiv 0$.

An example of complete smooth manifold $(M,g)$ with nonpositive Ricci curvature is the well know *Cartan-Hadamard manifold*, that is a simply connected complete Riemannian manifold of nonpositive sectional curvature. The following assertion is valid.

**Corollary 3.2**. *Let $(M,g)$ be a Cartan-Hadamard manifold of dimensional $n \geq 2$ with irreducible holonomy group $\text{Hol}(g)$. Then each infinitesimal harmonic transformation on $(M,g)$ with finite kinetic energy identically vanishes.*

**Remark 3.4**. Other properties of infinitesimal harmonic transformations can be found in our papers [20]; [5]; [13]; [15] and [16]. In particular, we proved in [20] that the set of all infinitesimal harmonic transformations on a compact Riemannian manifold $(M,g)$ is a finite-dimensional vector space over $\mathbb{R}$. Moreover, the Lie algebra of infinitesimal isometric transformation is a subspace of this vector space ([5]). We recall here that an infinitesimal isometric transformation (infinitesimal isometry) or a Killing vector field $X$ on $(M,g)$ is defined by the well known equation $L_X g = 0$.

The following theorem on infinitesimal isometric transformations is a well known oldest result (see, for example, [9, p. 57] and [19, p. 44]).

**Theorem 3.4**. *Let $(M,g)$ be a Riemannian manifold and $\xi$ a vector field on $(M,g)$. If $\xi$ is an infinitesimal isometry, it satisfies the following differential equations:*

$$(3.5) \qquad \tilde{\Delta}\theta = 2\,Ric(\xi,\cdot);$$

$$(3.6) \qquad div\,\xi = 0$$

*for the 1-form $\theta$ corresponding to $\xi$ under the duality defined by the metric $g$. Conversely, if $(M,g)$ is compact and $\xi$ satisfies (3.5) and (3.6), then $\xi$ is an infinitesimal isometry.*

The equation (3.6) is more rigid than required for the theorem above. In turn, we can formulate and prove an alternative version of Theorem 3.4.

**Theorem 3.5.** *Let $(M,g)$ be a Riemannian manifold and $\xi$ a vector field on $(M,g)$. If $\xi$ is an infinitesimal isometry, it satisfies the following conditions: $\xi$ is an infinitesimal harmonic transformation and*

(3.7) $$L_\xi \, div\,\xi \geq 0.$$

*Conversely, if $(M,g)$ is compact and $\xi$ satisfies the above two conditions, then $\xi$ is an infinitesimal isometry.*

**Proof.** We have already shown that if $\xi$ is an infinitesimal harmonic transformation, it satisfies (3.5). If, in addition, $\xi$ is an infinitesimal isometry, then the equation $div\,\xi = 0$ holds. This implies (3.7). To prove the converse, we may assume that $M$ is compact and orientable. (If $M$ is not orientable, consider its orientable double covering). Let's consider the vector field $X = (div\,\xi)\xi$ for an arbitrary infinitesimal harmonic transformation $\xi$ on $(M,g)$. The divergence of $X$ has the form

(3.8) $$div\,X = L_\xi\,(div\,\xi) + (div\,\xi)^2.$$

If we apply classic Green's theorem $\int_M div\,X \, dVol_g = 0$ to $X = (div\,\xi)\xi$, then we obtain the integral formula

(3.9) $$\int_M \left(L_\xi(div\,\xi) + (div\,\xi)^2\right) dVol_g = 0$$

for the canonical measure $dVol_g$ which is associated to the metric $g$. If the inequality $L_\xi(div\,\xi) \geq 0$ holds anywhere on $(M,g)$, then from (3.9) we conclude that $div\,\xi = 0$. Next, for complete the proof we can refer to Theorem 3.4.

**Remark 3.5.** The divergence of each smooth vector field $\xi$ on $(M,g)$ is a scalar function defined by (see, for example, [19, p. 6] and [21, p. 195])

(3.10) $$(div\,\xi)dVol_g = L_\xi(dVol_g).$$

Due to (3.10), the function $div\,\xi$ is called in [21, p. 195] the *logarithmic rate of volumetric expansion* along the flow generated by the vector field $\xi$. Therefore,

$L_\xi \, div\, \xi$ measures the *acceleration of volumetric expansion*, i.e. the acceleration of change of the volume element $dVol_g$ along trajectories of the flow with the velocity vector $\xi$. In particular, the condition $L_\xi \, div\, \xi \geq 0$ means that $dVol_g$ is a nondecreasing function along trajectories of this flow.

The equality $L_\xi \, g = 0$ also implies the invariance condition for the Ricci tensor $L_\xi \, Ric = 0$. In General Relativity, there are investigations (see, for example, [22] and [23]), where "weakened condition" of the form $trace_g(L_\xi \, Ric) = 0$ is studied instead of the condition $L_\xi \, Ric = 0$. It is of interest to note, that, in the case of compact $(M, g)$, the addition of conditions $L_\xi \, Ric \geq 0 \, (\leq 0)$ to equation (3.1) implies that infinitesimal harmonic transformation $\xi$ will actually be an infinitesimal isometry. Namely, the following theorem holds.

**Theorem 3.6**. *Let $(M, g)$ be a Riemannian manifold and $\xi$ a vector field on $(M, g)$. If $\xi$ is an infinitesimal isometry, it satisfies the following conditions: $\xi$ is an infinitesimal harmonic transformation and*

(3.11) $$trace_g \, L_\xi \, Ric \geq 0 \, (\leq 0).$$

*Conversely, if $(M, g)$ is compact and $\xi$ satisfies the above two conditions, then $\xi$ is an infinitesimal isometry.*

**Proof.** We have already shown that if $\xi$ is an infinitesimal harmonic transformation, it satisfies (3.5) that equals to the first condition of our theorem. If, in addition, $\xi$ is an infinitesimal isometry, then the equation $L_\xi \, Ric = 0$ holds. To prove the converse, we may assume that $M$ is compact. If $\xi$ is an infinitesimal harmonic transformation then we have the differential equation ([24])

$$\Delta \, div\, \xi = trace_g(L_\xi \, Ric).$$

Using the *Bochner maximum principal* ([9, Theorem 2.2]), we can conclude that $div\, \xi = const$. On the other hand, $\int_M div\, \xi \, dVol_g = 0$. Hence, $div\, \xi = 0$, showing that $\xi$ is an infinitesimal transformation (see Theorem 3.4). Theorem 3.6 is proved.

**Remark 3.6**. Theorem 3.6 is the second alternative version of the classic Theorem 3.4 on infinitesimal isometric transformations.

## 4. Ricci solitons

The main results of this section are applications to the Ricci soliton theory of the results of the third section of our paper.

Let $g$ be a fixed Riemannian metric on a smooth manifold $M$ and $Ric$ be its Ricci tensor. Consider the one-parameter family of diffeomorphisms $\varphi_t(x): M \to M$ that is generated by the smooth vector field $\xi$ on $M$. The evolutive metric $g(t) = \sigma(t)\varphi_t^*(x)g(0)$ for a positive scalar $\sigma(t)$ such that $\sigma(0) = 1$ and $g(0) = g$ is a *Ricci soliton* iff the metric $g$ is a solution of the nonlinear stationary PDF

$$(4.1) \qquad -2\,Ric = L_\xi\, g + 2\lambda\, g$$

where $L_\xi\, g$ is the *Lie derivative* of $g$ with respect to $\xi$ and $\lambda$ is a constant (see, for example, [4, p. 22]). To simplify the notation, we denote the Ricci soliton in the following way $(g, \xi, \lambda)$. A Ricci soliton is called *steady*, *shrinking* and *expanding* if $\lambda = 0$, $\lambda < 0$ and $\lambda > 0$, respectively. In addition, a Ricci soliton is called *Einstein* if $L_\xi\, g = 0$, and it is called *trivial* if $\xi \equiv 0$.

In [14], we have shown that the following theorem is true.

**Theorem 4.1**. *The vector field $\xi$ of an arbitrary Ricci soliton $(g, \xi, \lambda)$ on a smooth manifold $M$ is an infinitesimal harmonic transformation on the Riemannian manifold $(M, g)$.*

Therefore, from Theorem 3.1 and Theorem 4.1 we can obtain the following statement.

**Corollary 4.1**. *Let $(g, \xi, \lambda)$ be a Ricci soliton on a smooth manifold $M$. If the Ricci tensor $Ric$ of $g$ is quasi-negative in a connected, open domain $U \subset M$ and the energy density of the flow generated by the vector field $\xi$ has a local maximum in some point of $U$, then $(g, \xi, \lambda)$ is an expanding Ricci soliton.*

**Proof.** Let $(g,\xi,\lambda)$ be a Ricci soliton on a smooth manifold $M$, then its vector field $\xi$ is an infinitesimal harmonic transformation on the Riemannian manifold $(M,g)$. If, in addition, the Ricci tensor $Ric$ of $g$ is quasi-negative in a connected, open domain $U \subset M$ and the energy density of the flow $e(\xi)$ generated by $\xi$ has a local maximum in some point $x$ of $U$, then $\xi \equiv 0$ in any point $y \in U$, by our Theorem 3.1. In this case, from (4.1) we obtain the equality $Ric_y = -\lambda g_y$ and hence $\lambda > 0$.

**Remark 4.1.** From Corollary 3.1 we can conclude that a Ricci soliton $(g,\xi,\lambda)$ with the quasi-negative Ricci curvature of $g$ is trivial on a compact smooth manifold $M$. In turn, from Theorem 3.2 we obtain the following

**Corollary 4.2.** *Let $(g,\xi,\lambda)$ be a Ricci soliton on a compact smooth manifold $M$. If $L_\xi s \leq 0$ for the scalar curvature $s$ of $g$, then $(g,\xi,\lambda)$ is trivial.*

**Proof.** Consider a Ricci soliton $(g,\xi,\lambda)$ on a compact smooth manifold $M$. From (4.1) we obtain $div\,\xi = -(s+n\lambda)$ for the scalar curvature $s = \text{trace}_g\,Ric$, then the acceleration of volumetric expansion of the flow generated by the vector field $\xi$ of the Ricci soliton $(g,\xi,\lambda)$ has the form

(4.2) $$L_\xi(div\,\xi) = -L_\xi s.$$

If $L_\xi s \leq 0$ then from Theorem 3.5 we obtain that $\xi$ is an infinitesimal isometry, i.e. $L_\xi g = 0$. On the other hand, an arbitrary Ricci soliton $(g,\xi,\lambda)$ on a compact smooth manifold $M$ is a gradient soliton, i.e. $\theta = grad\,u$ for some $u \in C^\infty M$ ([7]). Therefore, the condition $L_\xi g = 0$ becomes $\nabla\nabla u = 0$ which implies the equation $\Delta u = 0$, i.e. $u$ is a *harmonic function*. In this case, $u = const$, by Bochner maximum principal. As a result, we obtain $\theta = grad\,u = 0$, and hence, the Ricci soliton $(g,\xi,\lambda)$ is trivial.

The following corollary of Theorem 3.6 is proved similarly.

**Corollary 4.3.** *Let $(g,\xi,\lambda)$ be a Ricci soliton on a compact smooth manifold $M$. If $\text{trace}_g\,L_\xi\,Ric \geq 0\;(\leq 0)$ for the Ricci tensor $Ric$ of $g$, then $(g,\xi,\lambda)$ is trivial.*

**Remark 4.2**. It is well known, that every steady and expanding Ricci soliton on a compact smooth manifold $M$ is trivial (see, for example, [7]). On the other hand, the well known problem presented in [7]: are there special conditions in dimensional $n \geq 4$ assuring that a shrinking compact Ricci soliton is trivial? Our two corollaries are answers to this problem.

The proof of the following corollary which is obtained from Theorems 3.3 is not required.

**Corollary 4.4**. *Let $(g,\xi,\lambda)$ be a Ricci soliton with complete Riemannian metric $g$ and nonpositive Ricci curvature on a connected smooth manifold $M$. If the kinetic energy of the flow generated by the vector field $\xi$ is infinite, then $(g,\xi,\lambda)$ is an Einstein Ricci soliton. If, in addition, the volume of $(M,g)$ is infinite or the Ricci curvature is negative at some point, then $(g,\xi,\lambda)$ is a trivial Ricci soliton.*

From (3.8) and (4.2) we conclude that the divergence of the acceleration vector of volumetric expansion $X = (div\,\xi)\xi$ has the form

(4.3) $$div\,X = -L_\xi s + (s + n\lambda)^2.$$

If, in addition, $(M,g)$ is a complete and oriented Riemannian manifold such that $\|(div\,\xi)\xi\| \in L^1(M,g)$ and $L_\xi s \leq 0$, then from (4.3) we obtain that $div\,X \geq 0$. Then, thanks to a *generalized Green's theorem* ([28]; [29]), we have $div\,X = 0$. It means that $L_\xi s = 0$ and

(4.4) $$s = -n\lambda.$$

In this case, from the well-known Schur´s identity $\delta\,Ric = -2^{-1}\nabla s$, we obtain the equation $\delta\,Ric = 0$. We recall the well-known equality $\delta\,\widetilde{\Delta} = \widetilde{\Delta}\,\delta$. Therefore, if we apply the divergence operator $\delta$ to both sides of the equation $\widetilde{\Delta}\theta = 2\,Ric(\xi,\cdot)$, we obtain $trace_g\,(S \cdot \nabla\xi) = 0$. In this case, from the equation (4.1) and the equality (4.4) we obtain the following $\|Ric\|^2 = n^{-1}s^2$. Then the square of the traceless part of the Ricci tensor is equal to zero, i.e. $\|Ric - n^{-1}s\,g\|^2 = \|Ric\|^2 - n^{-1}s^2 = 0$. Therefore, we conclude that $Ric = n^{-1} s \cdot g$ and $L_\xi g = 0$. This means that $\xi$ is an infinitesimal isometry and $(g, \xi, \lambda)$ is an Einstein soliton. Thus, we proved the following statement.

**Theorem 4.2**. *Let $(g, \xi, \lambda)$ be a Ricci soliton with complete Riemannian metric $g$ on a connected and oriented smooth manifold M such that*

*(i) $L_\xi s \leq 0$ for the scalar curvature s of the metric $g$;*

*(ii) $\|X\| \in L^1(M, g)$ for the logarithmic rate of volumetric expansion $X = (div\,\xi)\xi$, then the flow generated by $\xi$ consists of isometric transformations and $(g, \xi, \lambda)$ is an Einstein soliton.*

Let $(g, \xi, \lambda)$ be a shrinking Ricci soliton with complete metric $g$ on a connected and oriented smooth manifold $M$, such that $\xi$ satisfies (i) and (ii), then from Theorem 4.2 we obtain the following inequality $Ric = n^{-1} s \cdot g > 0$ for the positive constant $s$. Then it follows from Myers' diameter bound that $(M, g)$ must be compact (see, for example, [30, pp. 171; 254]). In particular, if $M$ is simply connected and $dim\,M = 3$, then $(M, g)$ be isometric a Euclidian sphere $\mathbb{S}^3$. In this case, we have the following

**Corollary 4.5**. *Let $(g, \xi, \lambda)$ be a shrinking Ricci soliton with complete Riemannian metric $g$ on a connected and oriented simply connected three-dimensional smooth manifold M such that $L_\xi s \leq 0$ for the scalar curvature s of the metric $g$ and $\|X\| \in L^1(M, g)$ for the logarithmic rate of volumetric expansion $X = (div\,\xi)\xi$, then $(M, g)$ be isometric a Euclidian sphere $\mathbb{S}^3$.*

In conclusion, Corollary 3.1 implies the following assertion.

**Corollary 4.6**. *Let $(g, \xi, \lambda)$ be a Ricci soliton with complete metric $g$, nonpositive Ricci tensor Ric and irreducible holonomy group $\mathrm{Hol}(g)$ on a connected noncompact smooth manifold M. If the kinetic energy of the flow generated by the vector field $\xi$ is infinite, then $(g, \xi, \lambda)$ is a trivial soliton.*

**Remark 4.3.** Our Theorem 4.2 generalizes Theorem 5.4 which is one of the main results of the paper [31]. In addition, our Corollary 4.4 is a generalization of Theorem 2 and Corollary 5 from our paper [14].

**Acknowledgments.** Our work was supported by Russian Foundation for Basis Research of the Russian Academy of Science (projects Nos. 16-01-00053 and 16-01-00756).

# References

[1] Eells, J., Sampson, J.H., Harmonic mappings of Riemannian manifolds, American Journal of Mathematics, **86** (1964), no. 1, 109-160.

[2] Schoen R., Yau S.T., Harmonic maps and the topology of stable hypersurfaces and smooth manifolds with nonnegative Ricci curvature, Comment. Math. Helvetici **61** (1976), no. 1, 333-341.

[3] Xin Y., Geometry of harmonic maps, Springer Science and Business Media, Boston and Berlin (2012).

[4] Nouhaud O., Transformations infinitésimales harmoniques, C. R. Acad. Sc. Paris, **274** (1972), 573–576.

[5] Stepanov S.E., Shandra I.G., Geometry of infinitesimal harmonic transformation, Ann. Glob. Anal. Geom., **24** (2003), 291-29.

[6] Chow B., Knopf D., The Ricci flow: An introduction, American Mathematical Society, Providence (2004).

[7] Eminenti M., La Nave G., Mantegazza C., Ricci solitons: The equation point of view, Manuscripta Mathematica, **127** (2008), no. 3, 345-367.

[8] Calabi E., An extension of E. Hopf's maximum principle with an application to Riemannian geometry, Duke Math. J., **25** (1957), 45-56.

[9] Bochner S., Yano K., Curvature and Betti numbers, Princeton, Princeton University Press (1953).

[10] Yau S.T., Some function-theoretic properties of complete Riemannian manifold and their applications to geometry, Indiana Univ. Math. J., **25** (1976), no. 7, 659-679.

[11] Li P., Schoen R., $L^p$ and mean value properties of subharmonic functions on Riemannian manifolds, Acta Mathematica, **153** (1984), 279-301.

[12] Stepanov S.E., Tsyganok I.I. Mikesh J., From infinitesimal harmonic transformations to Ricci solitons, Mathematica Bohemica, **1** (2013), 25-36.

[13] Stepanov S.E., Shandra I.G., Harmonic diffeomorphisms of smooth manifolds, St. Petersburg. Math. J., **16** (2005), no. 2, 401–412.

[14] Stepanov S.E., Shelepova V.N., A note on Ricci soliton, Math. Notes **86** (2009), no. 3, 447-450.

[15] Stepanov S.E., Tsyganok I.I., Infinitesimal harmonic transformations and Ricci solitons on complete Riemannian manifolds, Russian Math. **54** (2010), no 3, 84-87.

[16] Stepanov S.E., Tsyganok I.I., Harmonic transforms of complete Riemannian manifolds, Math. Notes, **100** (2016), no. 3, 465-471.

[17] Wu H., A remark on the Bochner technique in differential geometry, Proc. Amer. Math. Soc., **78** (1980), no. 3, 403-408.

[18] Stepanov S.E., Tsyganok I.I., Vanishing theorems for harmonic mappings into nonnegativity curved smooth manifolds and their applications, Manuscripta Math., **154** (2017), Issue 1-2, 79-90.

[19] Koboyashi S., Transformation groups in differential geometry, Springer-Verlag, Berlin and Heidelberg (1995).

[20] Stepanov S.E., Mikes J., The spectral theory of the Yano Laplacian with some of its applications, Ann. Glob. Anal. Geom. **48** (2015), 37-46.

[21] O´Neil B., Semi-Riemannian geometry with applications to Relativity, Academic Press, San Diego (1983).

[22] Davis W.R., Oliver D.R., Matter field space times admitting symmetry mappings satisfying vanishing contraction of the Lie deformation of the Ricci tensor, Ann. Inst. H. Poincare, Sec A (N.S.), **28** (1978), no. 2, 197-206.

[23] Green L.H., Norris L.K., Oliver D.R., Davis W.R., The Robertson-Walker metric and the symmetries belong to the family of contracted Ricci collineations, General Relativity and Gravitation, **8** (1997), no. 9, 731-736.

[24] Stepanov S.E., Shandra I.G., New characteristics of infinitesimal isometry and nontrivial Ricci solitons, Math. Notes **92** (2012), no. 3, 119-122.

[25] Arnold V.I., Kresin B.A., Topological methods in hydrodynamics, Springer-Verlag, New York (1998).

[26] Bérard P.H., From vanishing theorems to estimating theorems: the Bochner technique revisited, Bulletin of the American Math. Soc., **19** (1988) no. 2, 371-406.

[27] Yau S.T., Erratum: Some function-theoretic properties of complete Riemannian manifold and their applications to geometry, (1976), 659-670, Indiana Univ. Math. J., **31** (1982), no. 4, 607.

[28] Caminha A., Souza P., Camargo F., Complete foliations of space forms by hypersurfaces, Bull. Braz. Math. Soc., New Ser., **41** (2010), no. 3, 339–353.

[29] Caminha A., The geometry of closed conformal vector fields on Riemannian spaces, Bull. Braz. Math. Soc., New Ser., **42** (2011), no. 2, 277–300.

[30] Petersen P., Riemannian Geometry, Springer, NY (2006).

[31] Catino G., Mastrolia P., Monticelli D., Rigoli M., Analytical and geometric properties of generic Ricci Solitons, Trans. Am. Math. Soc. **368** (2016), no. 11, 7533–7549.

DEPARTMENT OF MATHEMATICS,
ALL RUSSIAN INSTITUTE FOR SCIENTIFIC
AND TECHNICAL INFORMATION OF THE
RUSSIAN ACADEMY OF SCIENCE,
125190 MOSCOW, RUSSIA
*E-mail address*: s.e.stepanov@mail.ru

DEPARTMENT OF MATHEMATICS, FINANCE UNIVERSITY
UNDER THE GAVERNMENT OF THE RUSSIAN FEDERATION,
125468 MOSCOW, RUSSIA
*E-mail address*: i.i.tsyganok@mail.ru